\documentclass[a4paper,10pt]{article}

\usepackage{amssymb}
\usepackage{amsmath}
\usepackage{amsthm}
\usepackage[dvips]{graphicx}
\usepackage{ascmac}
\usepackage{color}

\usepackage[dvipdfm,colorlinks=true,bookmarks=true,
bookmarksnumbered=true,bookmarkstype=toc,linktocpage=true
]{hyperref}

\newtheorem{thm}{Theorem}[section]

\newtheorem{ass}[thm]{Assumption}

\newcommand{\mca}{\mathcal{A}}\newcommand{\mcb}{\mathcal{B}}
\newcommand{\mcc}{\mathcal{C}}
\newcommand{\mcf}{\mathcal{F}}

\newcommand{\mcl}{\mathcal{L}}
\newcommand{\mcm}{\mathcal{M}}

\newcommand{\mbbn}{\mathbb{N}}

\newcommand{\mbbr}{\mathbb{R}}

\newcommand{\mbbrp}{\mathbb{R}_{+}}

\newcommand{\mbF}{\mathbf{F}}

 \newcommand{\lam}{\lambda} \newcommand{\ep}{\epsilon} 
\newcommand{\vp}{\varphi}

  \newcommand{\Sig}{\Sigma}
\newcommand{\Lam}{\Lambda}  

\newcommand{\p}{\partial}  

\def\nn{\nonumber}

\title{Edgeworth expansion for the integrated L\'evy driven Ornstein-Uhlenbeck process}
\date{Version: \today}
\author{Hiroki Masuda
\footnote{Corresponding author}
\footnote{
Institute of Mathematics for Industry, Kyushu University. 
Motooka 744, Nishi-ku, Fukuoka 819-0395, Japan.
Email: hiroki@imi.kyushu-u.ac.jp
} \and 
Nakahiro Yoshida
\footnote{
Graduate School of Mathematical Sciences, University of Tokyo. 
Komaba 3-8-1, Meguro-ku, Tokyo 153-8914, Japan.
}
}

\begin{document}
\maketitle

\begin{abstract}
We verify the Edgeworth expansion of any order for the integrated ergodic L\'evy driven Ornstein-Uhlenbeck process, 
applying a Malliavin calculus with truncation over the Wiener-Poisson space. 
Due to the special structure of the model, the coefficients of the expansion can be given in a closed form. 

\medskip

\end{abstract}


\section{Introduction}

Let $(X,Y)=\{(X_{t},Y_{t})\}_{t\in\mbbrp}$ be the bivariate model described by
\begin{equation}
\left\{
\begin{array}{l}
\displaystyle{X_{t}=X_{0}-\lam\int_{0}^{t}X_{s}ds+Z_{t},} \\
\displaystyle{Y_{t}=\int_{0}^{t}(\gamma+\beta X_{s})ds+\rho Z_{t},}
\end{array}
\right.
\label{basic-model}
\end{equation}
where $Z=(Z_{t})_{t\in\mbbrp}$ is a non-trivial L\'evy process independent of the initial variable $X_{0}$, 
and the parameter 
$(\lam,\gamma,\beta,\rho)\in(0,\infty)\times\mbbr\times(\mbbr\backslash\{0\})\times\mbbr$ satisfies that
\begin{equation}
\beta+\rho\lam\ne 0.
\label{para_cond}
\end{equation}
The process $X$ is the exponentially ergodic L\'evy driven Ornstein-Uhlenbeck (OU) process; 
we refer to \cite{Mas04} and the references therein for fundamental facts concerning the OU process. 
The goal of this note is to provide conditions under which 
the Edgeworth expansion of the expectation $E[f(T^{-1/2}H_{T})]$ as $T\to\infty$ is valid, where
\begin{equation}
H_{T}:=Y_{T}-E[Y_{T}] 
\label{H-def}
\end{equation}
and $f:\mbbr\to\mbbr$ is a measurable function of at most polynomial growth. 
The condition (\ref{para_cond}) will turn to be necessary for the Gaussian limit of $\mcl(T^{-1/2}H_{T})$ to be non-degenerate: 
as a matter of fact, the necessity of (\ref{para_cond}) can be seen concisely by the expression
\begin{equation}
T^{-1/2}H_{T}
=(\beta+\rho\lam)T^{-1/2}\int_{0}^{T}(X_{t}-E[X_{t}])dt
+\rho T^{-1/2}\left\{(X_{t}-E[X_{t}])
-(X_{0}-E[X_{0}])\right\},
\nonumber
\end{equation}
so that, if $\beta+\rho\lam=0$ and $(X_{t}-E[X_{t}])-(X_{0}-E[X_{0}])=O_{p}(1)$ as $T\to\infty$, 
then $\mcl(T^{-1/2}H_{T})$ tends in probability to $0$ (See Section \ref{sec_coeff}).

As is well known, distributional regularity of the underlying model is essential to 
the validity of the Edgeworth expansion. 
At first glance, the regularity of the joint distribution $\mcl(X,H)$, which will play an essential role 
in derivation of the expansion (see Section \ref{sec_proof}), does not seem enough 
since we have only one-dimensional random input $Z$ against the two-dimensional objective $(X,H)$. 
In particular, for pure-jump $Z$ we have to take distributional regularity 
over the Poisson space into account, rendering the problem mathematically interesting in its own right. 
In this case, we will execute the Malliavin calculus under truncation, 
which enables us to successfully pick out a nice event on which the integration by parts formula can apply to ensure 
distributional regularity; 
more specifically, our truncation functional will be constructed through two diffusive jumps, 
so as to make the Malliavin covariance matrix associated with the flow of $(X,H)$ non-degenerate 
(As will be mentioned in Section \ref{sec_truncation}, a single jump is not enough). 
The Malliavin calculus conveniently enables us to bypass intractable 
direct estimate of the characteristic function of $\mcl(T^{-1/2}H_{T})$, 
and results in fairly simple conditions.

Our result has the following statistical implication. 
Suppose that we can directly observe $\{X_{t}: 0\le t\le T\}$, 
based on which we want to estimate $\theta_{0}:=E[X_{0}]$ (the mean of the stationary distribution). 
A natural estimator is then given by
\begin{equation}
\hat{\theta}_{T}:=\frac{1}{T}\int_{0}^{T}X_{s}ds
\nonumber
\end{equation}
We easily see that 
$T^{-1/2}H_{T}=T^{1/2}(\hat{\theta}_{T}-\theta_{0})$ with $\beta=1$ and $\gamma=\rho=0$, 
hence the consistency, asymptotic normality, and 
higher order expansion of $\hat{\theta}_{T}$ are obtained according to our result.

\medskip

The main result is given in Section \ref{sec_results}, followed by the proof in Section \ref{sec_proof}.


\section{Edgeworth expansion}\label{sec_results}

\subsection{Statement of result}

We are given a stochastic basis $(\Omega,\mcf,\mbF=(\mcf_{t})_{t\in\mbbrp},P)$, 
on which our processes are defined.

\begin{ass}
$X$ is strictly stationary with a stationary distribution $F\in\bigcap_{p>0}L^{p}(P)$.
\label{ass1}
\end{ass}

We remark that: under Assumption \ref{ass1} $X$ is exponentially $\beta$-mixing and ergodic; 
Assumption \ref{ass1} is equivalent to $Z_{1}\in\bigcap_{p>0}L^{p}(P)$. 
See \cite{Mas04} for details.

Denote by $(b,C,\Pi)$ the generating triplet of $Z$ in the form
\begin{equation}
\vp(u;Z_{t})=
\exp\bigg\{t\bigg(ibu-\frac{1}{2}Cu^{2}+
\int_{\mbbr}(e^{iuz}-1-iuz)\Pi(dz)\bigg)\bigg\},
\nonumber
\end{equation}
where $b\in\mbbr$, $C\ge 0$, and the L\'evy measure $\Pi$ defined on $\mbbr$ 
is a $\sigma$-finite measure satisfying $\Pi(\{0\})=0$ and 
$\int_{0<|z|\le 1}z^{2}\Pi(dz)<\infty$. Then the process $H$ of (\ref{H-def}) satisfies
\begin{equation}
dH_{t}=\beta(X_{t}-\kappa_{F}^{(1)})dt+\rho d\bar{Z}_{t},\qquad H_{0}=0,
\nn
\end{equation}
where $\bar{Z}_{t}:=Z_{t}-E[Z_{t}]=Z_{t}-E[Z_{1}]t$ and
\begin{equation}
\kappa_{\xi}^{(k)}:=i^{-k}\p_{u}\log E[\exp(iu\xi)],
\nonumber
\end{equation}
the $k$-th cumulant of $\xi$, 
with $\p_{v}$ denoting the (partial) differentiation with respect to a variable $v$.

Denote by $\Lambda$ the Poisson random measure associated with jumps of $Z$. We decompose it as
\begin{equation}
\Lambda(dt,dz)=\mu^{\flat}(dt,dz)+\mu(dt,dz)
\nn
\end{equation}
for some Poisson random measures $\mu^{\flat}$ and $\mu$; 
by the independently scattered property of $\Lam$, such a decomposition is always possible. 
Correspondingly, we write
\begin{equation}
\Pi(dz)=\nu^{\flat}(dz)+\nu(dz),
\nn
\end{equation}
where $\nu^{\flat}$ and $\nu$ stand for the L\'evy measures on $\mbbrp$ 
associated with $\mu^{\flat}$ and $\mu$, respectively.

\begin{ass}
Either one of the following two conditions holds true:
\begin{itemize}
\item[{\rm (i)}] $C>0$ (no condition is imposed on the jump-part characteristic);
\item[{\rm (ii)}] $C=0$ and there exists a non-empty open subset of 
$\mbbr\backslash\{0\}$ on which $\nu$ admits a positive $\mcc^{3}$-density, say $g$,  
with respect to the Lebesgue measure.
\end{itemize}
\label{ass2}
\end{ass}

Let us introduce the notation necessary for the Edgeworth expansion; see \cite{Yos04} for more details. 
We introduce the $r$-th cumulant function of $T^{-1/2}H_{T}$ ($r\in\mbbn$, $r\ge 2$):
\begin{equation}
\chi_{r,T}(u):=\partial_{u}^{r}\log E\left[\exp(iuT^{-1/2}H_{T})\right].
\nonumber
\end{equation} 
Let $p\ge 3$ be an integer. 
The {\it $(p-2)$-th Edgeworth expansion $\Psi_{p,T}$} (a signed measure) 
is defined by the Fourier inversion of $u\mapsto\hat{\Psi}_{p,T}(u)$, 
where 
\[\hat{\Psi}_{p,T}(u):=\exp\left(\frac{1}{2}\chi_{T,2}(u)\right)
+\sum_{r=1}^{p-2}T^{-r/2}\tilde{P}_{r,T}(u),\]
with $\tilde{P}_{r,T}(u)$ specified via the formal expansion
\begin{equation*}
\exp\bigg(\sum_{r=2}^{\infty}\frac{1}{r!}\chi_{r,T}(u)\bigg)
=\exp\left(\frac{1}{2}\chi_{2,T}(u)\right)
+\sum_{r=1}^{\infty}T^{-r/2}\tilde{P}_{r,T}(u).
\end{equation*}
Let $\phi(\cdot;\Sigma)$ stand for the one-dimensional centered normal density 
having variance $\Sig>0$, then 
the $r$-th Hermite polynomial associated with $\phi(\cdot;\Sigma)$ is 
$h_{r}(y;\Sigma):=(-1)^{r}\phi(y;\Sigma)^{-1}\p_{y}^{r}\phi(y;\Sigma)$. 
Let 
\begin{equation}
\chi_{r,T}:=(-i)^{r}\chi_{r,T}(0),
\nonumber
\end{equation}
the $r$-th cumulant of $T^{-1/2}H_{T}$; in Section \ref{sec_coeff}, 
we will see that $\chi_{r,T}=O(T^{-(r-2)/2})$ as $T\to\infty$. 
The density of $\Psi_{p,T}$ with respect to the Lebesgue measure is given by
\begin{equation*}
g_{p}(y;T^{-1/2}H_{T})=\left\{1+
\sum_{k=1}^{p-2}\sum_{l=1}^{k}\sum_{k_{1},\dots,k_{l}\in\mbbn:\atop 
k_{1}+\cdots +k_{l}=k}\frac{\chi_{k_{1}+2,T}\cdot\cdots\cdot\chi_{k_{l}+2,T}}
{l!(k_{1}+2)!\cdots (k_{l}+2)!}h_{k+2l}(y;\Sigma_{T})
\right\}\phi(y;\Sigma_{T}),
\end{equation*}
where $\Sigma_{T}:=\chi_{2,T}$; we will approximate $E[f(T^{-1/2}H_{T})]$ by 
$\Psi_{p,T}[f]:=\int f(y)g_{p}(y;T^{-1/2}H_{T})dy$. 

Let $p_{0}:=2[p/2]$ and denote by $\mathcal{E}(M,p_{0})$ the set of all measurable 
functions $f:\mbbr\to\mbbr$ satisfying $|f(x)|\le M(1+|x|^{p_{0}})$ for every $x\in\mbbr$.

Now we can state the main result. 

\begin{thm}
Let $X,Y,H$ be given through (\ref{basic-model}) and (\ref{H-def}), 
and suppose that (\ref{para_cond}) and Assumptions \ref{ass1} and \ref{ass2} hold true. 
Fix any positive number $\Sigma^{0}$ such that
\[\Sigma^{0}>\frac{2}{\lam}(\beta+\rho\lam)^{2}\kappa_{F}^{(2)}.\] 
Then, for any $M,K>0$, there exist positive constants $M^{\ast}$ and 
$\delta^{\ast}$ such that
\begin{equation}
\left|E[f(T^{-1/2}H_{T})]-\Psi_{p,T}[f]\right|
\le M^{\ast}
\int_{\mbbr}\sup_{|y|\le T^{-K}}|f(x+y)-f(x)|\phi(x;\Sigma^{0})dx
+o(T^{-(p-2+\delta^{\ast})/2})
\label{ae}
\end{equation}
for $T\to\infty$ uniformly in $f\in\mathcal{E}(M,p_{0})$.
\label{th2}
\end{thm}

Most often in practice, the first term in the upper bound in (\ref{ae}) can be quickly vanishing by taking $K$ large; 
for example, it is the case when $f$ is an indicator function $f=1_{A}$ for various $A\subset\mbbr$, 
such as $A=(-\infty,a]$, $A=[a,b]$, and so on.


\subsection{Explicit coefficients}\label{sec_coeff}

The approximating density $g_{p}(\cdot;T^{-1/2}H_{T})$ involves the cumulants $\chi_{2,T},\chi_{3,T},\dots,\chi_{p,T}$. 
We here prove the explicit formula for $\chi_{r,T}$, $r\ge 2$.

Noticing the explicit solution $X_{t}=e^{-\lam t}X_{0}+\int_{0}^{t}e^{-\lam (t-s)}dZ_{s}$, 
we can apply the stochastic Fubini theorem to obtain the relation
\begin{equation}
\int_{0}^{t}X_{s}ds=\eta(\lam,t)X_{0}+\int_{0}^{t}\eta(\lam,t-s)dZ_{s},
\label{intOU-SIE}
\end{equation}
where $\eta(\lam,u)=\lam^{-1}(1-e^{-\lam u})$; 
one can consults \cite{BS03} for a detailed analysis of integrated OU processes, 
especially in the context of financial econometrics. 
It follows from (\ref{basic-model}), (\ref{intOU-SIE}), 
and the special relation $k\lam\kappa_{F}^{(k)}=\kappa_{Z_{1}}^{(k)}$ for $k\in\mbbn$ (see \cite{bs01a,Mas04}) 
that we can express $H_{T}$ as
\begin{equation}
H_{T}=\beta\eta(\lam,T)X_{0}-T(\beta+\rho\lam)\kappa^{(1)}_{F}
+\int_{0}^{T}\left\{\rho+\beta\eta(\lam,T-s)\right\}dZ_{s}.
\nonumber
\end{equation}
Hence, using the independence between $X_{0}$ and $Z$ we obtain
\begin{align}
\chi_{r,T}&=
(-i)^{r}\left[\partial_{u}^{r}\kappa\left(\beta T^{-1/2}\eta(\lam,T)u;F\right)
+\int_{0}^{T}\partial_{u}^{r}\kappa
\left(\{\rho+\beta\eta(\lam,T-s)\}T^{-1/2}u;Z_{1}\right)ds\right]\bigg|_{u=0}
\nonumber\\
&=\left\{\beta T^{-1/2}\eta(\lam,T)\right\}^{r}\kappa^{(r)}_{F}
+\int_{0}^{T}\left(\{\rho+\beta\eta(\lam,v)\}T^{-1/2}\right)^{r}dv\lam r\kappa^{(r)}_{F}
\nn\\
&=T^{-(r-2)/2}\left[T^{-1}\left\{\beta \eta(\lam,T)\right\}^{r}
+\lam rT^{-1}\int_{0}^{T}\left\{\rho+\beta\eta(\lam,v)\right\}^{r}dv
\right]\kappa^{(r)}_{F}.
\nn
\end{align}
By making use of the differential equation 
$\p_{s}\left\{\eta(\lam,s)\right\}^{k}=k\left\{\eta(\lam,s)\right\}^{k-1}
-\lam k\left\{\eta(\lam,s)\right\}^{k}$ with $\eta(\lam,0)=0$ and then 
integrating the both sides with respect to $s$ over $[0,T]$, 
we can proceed as in \cite[Section 3]{MasYos05} to conclude that
\begin{equation}
\chi_{r,T}=T^{-(r-2)/2}\left[T^{-1}\left\{\beta \eta(\lam,T)\right\}^{r}
+\lam r\sum_{j=0}^{r}\binom{r}{j}\rho^{r-j}\beta^{j}\mcm_{r,T}(j)
\right]\kappa^{(r)}_{F},
\label{eq+2}
\end{equation}
where $\mcm_{r,T}(j)$ is given by
\begin{align}
\mcm_{r,T}(0)&=1, \nn\\
\mcm_{r,T}(j)&=\lam^{-j}-T^{-1}\lam^{-(j+1)}\sum_{k=1}^{j}k^{-1}\left\{\lam\eta(\lam,T)\right\}^{k},\quad j\ge 1.
\nonumber
\end{align}
Thus we can explicitly write down the coefficients of the Edgeworth expansion $\Psi_{p,T}$ up to any order. 
It is obvious from (\ref{eq+2}) that $\chi_{r,T}=O(T^{-(r-2)/2})$ for $r\ge 2$;
\begin{equation}
T^{(r-2)/2}\chi_{r,T}\to
\lam r\sum_{j=0}^{r}\binom{r}{j}\rho^{r-j}\beta^{j}\lam^{-j}
\kappa^{(r)}_{F}.
\nonumber
\end{equation}
In particular,
\begin{equation}
\Sigma_{T}=\chi_{2,r}\to 2\lam^{-1}(\beta+\rho\lam)^{2}\kappa_{F}^{(2)},
\nn
\end{equation}
hence the necessity of the condition (\ref{para_cond}).


\section{Proof of Theorem \ref{th2}}\label{sec_proof}

We will apply \cite[Theorem 1]{Yos04}. 
In order to ensure distributional regularity necessary for the Edgeworth expansion, 
we will make use of a Malliavin calculus with an effective truncation functional. 
The main idea of the proof is in principle similar to that of \cite[Section 4]{MasYos05} 
treating the stochastic volatility model, where $X$ expresses the latent positive volatility process. 
However, the OU process $X$ in the present model can take negative values too, and, as such, 
the way of constructing a truncation functional is essentially different from 
that of \cite{MasYos05}. 
To save space, we will sometimes omit the technical details, referring to the pertinent parts of \cite{bgj,MasYos05}.

Let us briefly overview the fundamental device. 
By means of \cite[Theorem 1]{Yos04}, it suffices to verify the following conditions:
\begin{itemize}
\item[{\rm $[A1]$}] {\it $X$ is strongly mixing with exponential rate};
\item[{\rm $[A2]$}] {\it $\sup_{t\in[0,T]}\|H_{t}\|_{L^{p+1}(P)}<\infty$ for each $T\in\mbbrp$};
\item[{\rm $[A3]$}] {\it there exist positive constants $t^{0}$, $a$, $a'$ and $B$, and a truncation functional 
$\psi: (\Omega,\mcf)\to\left([0,1],\mcb([0,1])\right)$ such that $0<a,a'<1$, $4a'<(a-1)^{2}$ and that}
\begin{eqnarray}
& &E\bigg[\sup_{|u|\ge B}\big|E[\psi\exp(iuH_{t^{0}})|X_{0},X_{t^{0}}]\big|\bigg]
<a', 
\nn
\\
& &1-E[\psi]<a. 
\nn
\end{eqnarray}
\end{itemize}
As was mentioned in Section \ref{sec_results}, Assumption \ref{ass1} ensures $[A1]$ and $[A2]$ (see (\ref{eq+2})), 
so that it remains to verify $[A3]$, which is a version of {\it conditional Cram\'er conditions}. 
Although it may be difficult in general to verify $[A3]$, 
we will be able to construct a specific truncation $\psi$ which significantly simplify the task. 

We also note that the condition $(\tilde{A}'-4)$ of \cite[p. 60 and p.130]{bgj} 
(smoothness of the coefficients, and integrability under cut-off through an auxiliary function) is indispensable. 
We will mention this point in Section \ref{sec_mcm}


\subsection{Transformation of the Poisson random measure}\label{sec_tprm}

In order to execute a Malliavin calculus of \cite{bgj}, we introduce a transformation of 
the absolutely continuous part of the Poisson random measure.

Under Assumption \ref{ass1}, $Z$ admits a L\'evy-It\^o decomposition of the form
\begin{equation*}
Z_{t}=\lam\kappa_{F}^{(1)}t+\sqrt{C}\tilde{w}_{t}
+\int_{0}^{t}\!\!\!\int_{\mbbr}z\tilde{\mu}^{\flat}(ds,dz)
+\int_{0}^{t}\!\!\!\int_{\mbbr}z\tilde{\mu}(ds,dz),\quad t\in\mbbrp,
\end{equation*}
where $\tilde{w}$ stands for a one-dimensional Wiener process defined on $(\Omega,\mcf,\mbF,P)$, 
$\tilde{\mu}^{\flat}(dt,dz):=\mu^{\flat}(dt,dz)-\nu^{\flat}(dz)dt$, 
and $\tilde{\mu}(dt,dz):=\mu(dt,dz)-\nu(dz)dt$.

Assumption \ref{ass2} assures the existence of a bounded domain
\begin{equation}
E_{0}=(c_{1},c_{2})\subset\mbbr\backslash\{0\},
\nonumber
\end{equation}
for which the L\'evy density $g$ of $\nu$ satisfies that
\begin{equation}
\inf_{z\in E_{0}}g(z)>0.
\nonumber
\end{equation}
Without loss of generality, we may and do suppose that 
$c_{1},c_{2}>0$: if $\nu(\mbbrp)\equiv 0$, then take $-Z$ as $Z$ anew. 
We introduce the change of variables $z^{\ast}=z^{\ast}(z)=g^{+}(z)$ 
through $z^{\ast}=z^{\ast}(z)=\int_{z}^{c_{2}}g(v)dv$ for $z\in E_{0}$; 
obviously, $g^{+}$ is strictly decreasing on $E_{0}$. 
Let $g^{-}$ denote the strictly decreasing inverse function of $g^{+}$ defined on
\begin{equation}
E=(g^{+}(c_{2}),g^{+}(c_{1})).
\nonumber
\end{equation}
Let $\mu^{\ast}$ denote the integer-valued random measure defined by
\begin{equation}
\int_{0}^{t}\!\!\!\int_{a_{1}}^{a_{2}}h(s,z)\mu(ds,dz)=
\int_{0}^{t}\!\!\!\int_{g^{+}(a_{2})}^{g^{+}(a_{1})}h(s,g^{-}(z^{\ast}))
\mu^{\ast}(ds,dz^{\ast})
\nn
\end{equation}
for each $t\in\mbbrp$, $a_{1},a_{2}\in\mbbr$ such that $a_{1}<a_{2}$, and 
for any measurable function $h$ on $\mbbrp\times\mbbrp$; in particular,
\begin{equation}
E[\mu^{\ast}([0,t],B)]=t{\rm Leb}(B).
\nonumber
\end{equation}
Writing $\tilde{\mu}^{\ast}(dt,dz^{\ast})=\mu^{\ast}(dt,dz^{\ast})-dtdz^{\ast}$, 
we transform $\mu$ (on $[0,t]\times E_{0}$) into $\mu^{\ast}$ as follows:
\begin{equation}
\int_{0}^{t}\!\!\!\int_{c_{1}}^{c_{2}}z\tilde{\mu}(ds,dz)=
\int_{0}^{t}\!\!\!\int_{g^{+}(c_{2})}^{g^{+}(c_{1})}g^{-}(z^{\ast})
\tilde{\mu}^{\ast}(ds,dz^{\ast}).
\nn
\end{equation}
The bivariate process $(X,H)$ satisfies the stochastic differential equation
\begin{equation}
\begin{split}
\binom{dX_{t}}{dH_{t}}&=(\kappa^{(1)}_{F}-X_{t})\binom{\lam}{-\beta}dt+
\sqrt{C}\binom{1}{\rho}d\tilde{w}_{t} \\
& {}\qquad+\int_{\mbbr}z\binom{1}{\rho}(\tilde{\mu}^{\flat}+1_{E_{0}^{c}}\tilde{\mu})(dt,dz)+
\int_{E\cup[g^{+}(c_{1}),\infty)}J(z^{\ast})\binom{1}{\rho}\tilde{\mu}^{\ast}(dt,dz^{\ast}),
\end{split}\label{sde-XH-B}
\end{equation}
where $J(z^{\ast}):=g^{-}(z^{\ast})1_{E}(z^{\ast})$ for 
$z^{\ast}\in E\cup[g^{+}(c_{1}),\infty)$. 
As $g^{-}$ is strictly decreasing, we have $|\p J(z^{\ast})|>0$ for $z^{\ast}\in E\cup[g^{+}(c_{1}),\infty)$.


\subsection{Malliavin covariance matrix}\label{sec_mcm}

Fix any constant $t^{0}>0$ and define $(\hat{\Omega},\hat{\mcb},\hat{P})$ 
to be the Wiener-Poisson canonical space (see \cite[the last paragraph in page 1178]{MasYos05}), 
on which we are given the flow $(X(\cdot,v),H(\cdot,v))^{\top}$ 
associated with $(X,H)$ of (\ref{sde-XH-B}) starting from $v=(x,h)^{\top}\in\mbbr^{2}$:
\begin{equation}
\begin{split}
\binom{X(t,v)}{H(t,v)}&=\binom{x}{h}
+\int_{0}^{t}(\kappa^{(1)}_{F}-X(s,v))\binom{\lam}{-\beta}ds+
\sqrt{C}\binom{1}{\rho}\tilde{w}_{t} \\
& {}\qquad+\int_{0}^{t}\int_{\mbbr}z\binom{1}{\rho}(\tilde{\mu}^{\flat}+1_{E_{0}^{c}}\tilde{\mu})(ds,dz)
+\int_{0}^{t}\int_{E\cup[g^{+}(c_{1}),\infty)}J(z^{\ast})\binom{1}{\rho}\tilde{\mu}^{\ast}(ds,dz^{\ast}).
\end{split}\nn
\end{equation}
Under the present assumption, the flow $(X(\cdot,\hat{v}),H(\cdot,\hat{v}))^{\top}$ 
clearly satisfies the condition $(\tilde{A}'-4)$.

Let $\hat{x}$ be a random variable independent of $(\tilde{w},\mu^{\flat}+1_{E_{0}^{c}}\mu,\mu^{\ast})$ 
such that $\mcl(\hat{x}|\hat{P})=F$ (the distribution under $\hat{P}$), and $\hat{v}:=(\hat{x},0)^{\top}$. 
We will compute the Malliavin covariance matrix of $(X(t^{0},\hat{v}),H(t^{0},\hat{v}))^{\top}$, 
whose ``non-degeneracy'' is essential here.

Let $Q\in\mbbr^{2}\otimes\mbbr^{2}$ be given by
\begin{equation*}
Q=\begin{pmatrix}
-\lam & 0 \\ \beta & 0
\end{pmatrix}.
\end{equation*}
In view of (\ref{sde-XH-B}), the process $K(t,v):=\p_{v}(X(t,v),H(t,v))^{\top}$ satisfies that, for each $v$,
\begin{align}
\frac{d}{dt}K(t,v)
&=\left(
\begin{array}{cc}
-\lam\p_{x}X(t,v) & 0 \\
\beta\p_{x}X(t,v) & 0
\end{array}\right)
=QK(t,v),
\nonumber
\end{align}
so that
\begin{equation*}
K(t^{0},\hat{v})=\exp(t^{0}Q)=
\begin{pmatrix}
e^{-\lam t^{0}} & 0 \\ \beta\lam^{-1}(1-e^{-\lam t^{0}}) & 1
\end{pmatrix}.
\end{equation*}
(Different from \cite[Eq.(29) in page 1181]{MasYos05}, $S(\cdot,\hat{v})$ is independent of $\hat{v}$.) 

Pick positive constants $c_{j}'$ and $c_{j}''$ ($j=1,2$) in such a way that 
$0<c_{1}<c_{1}'<c_{1}''<c_{2}''<c_{2}'<c_{2}<\infty$, and let
\begin{equation}
\check{E}:=(g^{+}(c_{2}''),g^{+}(c_{1}'')).
\nonumber
\end{equation}
(Trivially, $\check{E}\Subset E$.) 
Let $\eta:\mbbrp\to\mbbrp$ be any bounded smooth function satisfying the conditions:
\begin{itemize}
\item[(i)] $\inf_{z^{\ast}\in\check{E}}\eta(z^{\ast})>0$;
\item[(ii)] $\eta(z^{\ast})=0$ for $z^{\ast}\notin(g^{+}(c_{2}'),g^{+}(c_{1}'))$.
\end{itemize}
The Malliavin covariance matrix of $(X(t^{0},\hat{v}),H(t^{0},\hat{v}))^{\top}$ is then well-defined and given by
\begin{equation}
U(t^{0},\hat{v}):=\exp(t^{0}Q)S(t^{0},\hat{v})\exp(t^{0}Q^{\top}),
\nonumber
\end{equation}
where
\begin{align}
S(t,\hat{v})&=C\int_{0}^{t}\exp(-sQ)
\begin{pmatrix}
1 & \rho \\ \rho & \rho^{2}
\end{pmatrix}
\exp(-sQ^{\top})ds
\nn\\
& {}\qquad+\int_{0}^{t}\!\!\!\int_{E}\exp(-sQ)
\begin{pmatrix}
1 & \rho \\ \rho & \rho^{2}
\end{pmatrix}
\exp(-sQ^{\top})V(z^{\ast})\mu^{\ast}(ds,dz^{\ast}),
\label{S_{t}-B}
\end{align}
with $V(z^{\ast}):=\{\p J(z^{\ast})\}^{2}\eta(z^{\ast})$; 
see \cite[Section 10]{bgj} for details of (\ref{S_{t}-B}). Thus we arrive at the identity
\begin{equation}
\textrm{det}U(t^{0},\hat{v})=e^{-2\lam t^{0}}\textrm{det}S(t^{0},\hat{v}),\quad\text{a.s.}
\label{U-B}
\end{equation}


\subsection{Completion of the proof under Assumption \ref{ass2} (i)}\label{sec_pC}

Suppose that $C>0$. It follows from (\ref{S_{t}-B}) that, in the matrix sense,
\begin{align}
S(t^{0},\hat{v})&\ge C\int_{0}^{t^{0}}e^{-sQ}
\begin{pmatrix}
1 & \rho \\ \rho & \rho^{2}
\end{pmatrix}
e^{-sQ^{\top}}ds
\nn\\
&=
\begin{pmatrix}
H_{2} & \text{sym.} \\
\chi H_{1}-(\beta/\lam)H_{2} & \chi^{2}t^{0}-2(\beta/\lam)\chi H_{1}+(\beta/\lam)^{2}H_{2}
\end{pmatrix}
\nonumber
\end{align}
where $H_{k}:=\int_{0}^{t^{0}}e^{k\lam s}ds$ and $\chi:=\rho+\beta/\lam$. 
The determinant of the rightmost side is
\begin{equation}
C^{2}\lam^{-4}(\beta+\rho\lam)^{2}\bigg\{\frac{\lam t^{0}}{2}
(e^{2\lam t^{0}}-1)-(e^{\lam t^{0}}-1)^{2}\bigg\},
\nn
\end{equation}
which is positive as soon as $t^{0}\lam\ne 0$ and $\beta+\rho\lam\ne 0$. 
Thus $S(t^{0},\hat{v})$ is bounded from below by a positive-definite matrix, hence the non-degeneracy of 
$U(t^{0},\hat{v})$ follows from (\ref{U-B}) without any non-trivial truncation functional; 
simply let $\psi\equiv 1$ in $[A3]$. 
Thus we have obtained the non-degeneracy of the Malliavin covariance matrix 
(i.e. enough integrability of $\{{\rm det}U(t^{0},\hat{v})\}^{-1}$), 
which corresponds to \cite[Lemma 6]{MasYos05}.

We further notice the following.
\begin{itemize}
\item The flow $(X(t,\hat{v}),H(t,\hat{v}))_{t\in[0,t^{0}]}$ satisfies the condition $(\tilde{A}'-4)$ 
(as was seen in Section \ref{sec_mcm}), 
hence the analogous assertions as \cite[Lemmas 7]{MasYos05} holds true.

\item Following the same argument as in \cite[pp.1184--1185]{MasYos05}, 
we see that there exists a random variable $\Phi'_{t^{0}}\in L^{1}(\hat{P})$ 
such that
\begin{equation}
E\left[\sup_{|u|\ge B}\left|E[\exp(iuH_{t^{0}})|X_{0},X_{t^{0}}]\right|\right]\le\frac{1}{B}\hat{E}[|\Phi'_{t^{0}}|]
\nonumber
\end{equation} 
for every $B>0$.
\end{itemize}

After all, we have deduced the analogous assertions to \cite[Lemmas 6, 7 and 8]{MasYos05}, 
completing the proof of Theorem \ref{th2} under Assumption \ref{ass1} and Assumption \ref{ass2} (i).


\subsection{Construction of a truncation functional}\label{sec_truncation}

It remains to prove Theorem \ref{th2} under Assumptions \ref{ass1} and \ref{ass2} (ii); 
then, in order to verify distributional regularity we have to make an effective use of jumps. 
We will construct the truncation functional $\psi$ in an explicit way through two diffusive jumps.

We continue the argument of Section \ref{sec_mcm}. 
Let $t_{1},t_{2}\in(0,t^{0})$ be constants such that $t_{1}<t_{2}$, and fix $z_{0}\in\check{E}$. 
Let $\ep>0$ be sufficiently small so that:
\begin{itemize}
\item $\overline{I^{\ep}_{1}}\cap\overline{I^{\ep}_{2}}=\emptyset$ 
for $I_{j}^{\ep}:=(t_{j}-\ep,t_{j}+\ep)$, $j=1,2$;
\item $g^{+}(c_{2}'')<z_{0}-\ep<z_{0}+\ep<g^{+}(c_{1}'')$.
\end{itemize}
Let $E^{\ep}:=(z_{0}-\ep,z_{0}+\ep)$ and
\begin{equation}
\mca^{\ep}:=\{\mu^{\ast}(I_{j}^{\ep},E^{\ep})=1\ 
\textrm{for}\ j=1,2.\}.
\label{setA-B}
\end{equation}
According to the independently scattered property of $\mu^{\ast}$ and 
since the L\'evy measure associated with $\mu^{\ast}$ over $E$) here is the Lebesgue one, we have
\begin{equation}
\hat{P}[\mca^{\ep}]=\left\{\hat{P}\left[\mu^{\ast}\left([0,2\ep],[0,2\ep]\right)\right]\right\}^{2}
=\left\{4\ep^{2}\exp(-4\ep^{2})\right\}^{2}>0
\nonumber
\end{equation}
for each $\ep>0$.

We define the truncation functional $\hat{\psi}_{\ep}$ by 
$\hat{\psi}_{\ep}=\zeta(\hat{\xi}_{\ep})$, 
where $\zeta:\mbbrp\to[0,1]$ is a non-increasing smooth function such that
$\zeta(x)=1$ if $0\le x\le 1/2$ and $\zeta(x)=0$ if $x\ge 1$, where
\begin{equation}
\hat{\xi}_{\ep}=\frac{2}{1+3\textrm{det}U(t^{0},\hat{v})}.
\label{xi-B}
\end{equation}
We will show that the Malliavin covariance matrix $U(t^{0},\hat{v})$ is 
non-degenerate on the event $\mca^{\ep}$ for any $\ep>0$ small enough.

Noting that
\begin{equation*}
\sup_{s: |s-t_{j}|\le\ep}\left|e^{-sQ}
-
\begin{pmatrix}
e^{\lam t_{j}} & 0 \\ \beta\lam^{-1}(1-e^{\lam t_{j}}) & 1
\end{pmatrix}
\right|
+
\sup_{z: |z-z_{0}|\le\ep}\left|
V(z)-V(z_{0})
\right|
\to 0
\end{equation*}
as $\ep\to 0$ and by virtue of (\ref{setA-B}), 
we apply Taylor's expansion around $z_{0}$ and $t_{j}$ ($j=1,2$) on $\mca^{\ep}$ to conclude that
\begin{align}
S(t^{0},\hat{v})
&\ge\sum_{j=1}^{2}\int_{I_{j}^{\ep}}\!\int_{E^{\ep}}e^{-sQ}
\begin{pmatrix}
1 & \rho \\ \rho & \rho^{2}
\end{pmatrix}
e^{-sQ^{\top}}V(z^{\ast})\mu^{\ast}(ds,dz^{\ast})
\nn\\
&=\sum_{j=1}^{2}e^{-t_{j}Q}
\begin{pmatrix}
1 & \rho \\ \rho & \rho^{2}
\end{pmatrix}
e^{-t_{j}Q^{\top}}V(z_{0})+o(1)
\nn\\
&=V(z_{0})M^{\ep}+o(1)
\nonumber
\end{align}
as $\ep\to 0$ (we use the symbol $o(1)$ for matrices too), where
\begin{equation}
M^{\ep}:=\begin{pmatrix}
J^{(2)} & \textrm{sym.} \\
(\rho+\beta\lam^{-1})J^{(1)}-\beta\lam^{-1}J^{(2)} & 
2(\rho+\beta\lam^{-1})^{2}-2\beta\lam^{-1}(\rho+\beta\lam^{-1})J^{(1)}+\beta^{2}\lam^{-2}J^{(2)}
\end{pmatrix}
\nonumber
\end{equation}
with $J^{(1)}:=e^{\lam t_{1}}+e^{\lam t_{2}}$ and $J^{(2)}:=e^{2\lam t_{1}}+e^{2\lam t_{2}}$. 
Therefore
\begin{equation*}
{\rm det}S(t^{0},\hat{v})\ge V(z_{0})^{2}\lam^{-2}(\beta+\lam\rho)^{2}
(e^{\lam t_{1}}-e^{\lam t_{2}})^{2}+o(1),
\end{equation*}
which is positive for $\ep$ sufficiently small whenever $\rho\lam+\beta\ne 0$ and $t_{1}\ne t_{2}$. 
[We note that a single jump is not enough: if we instead estimate $S(t^{0},\hat{v})$ as
\begin{equation}
S(t^{0},\hat{v})\ge e^{-t_{1}Q}
\begin{pmatrix}
1 & \rho \\ \rho & \rho^{2}
\end{pmatrix}
e^{-t_{1}Q^{\top}}V(z_{0})+o(1),
\nonumber
\end{equation}
then the determinant of the first term in the right-hand side turns out to be identically $0$.]

We may set $V(z_{0})$ arbitrarily large by choosing the function $\eta$ suitably, so that, 
recalling (\ref{U-B}) we conclude that
$\textrm{det}U(t^{0},\hat{v})\ge 1$ on $\mca^{\ep}$ for some $\ep>0$. 
Fix such such $\eta$ and $\ep>0$. The definition (\ref{xi-B}) leads to the estimate
\begin{equation}
\hat{P}[\hat{\xi}_{\ep}\le 1/2]
\ge\hat{P}\left[\left\{\textrm{det}U(t^{0},\hat{v})\ge 1\right\}\cap\mca^{\ep}\right]
=\hat{P}[\mca^{\ep}]>0,
\nonumber
\end{equation}
hence the assertion corresponding to \cite[Lemma 6]{MasYos05} holds true. 
We fix the $\eta$ and $\ep>0$ in the sequel.

Clearly $\hat{\psi}_{\ep}>0$ implies that $1/3\le\textrm{det}U(t^{0},\hat{v})$, hence
\begin{equation}
\hat{\psi}_{\ep}\left\{{\rm det}U(t^{0},\hat{v})\right\}^{-1}
\in\bigcap_{0<p<\infty}L^{p}(\hat{P}).
\nonumber
\end{equation}
This implies that the integration-by-parts formula under the truncation $\hat{\psi}_{\ep}$ is in force. 
Then, as before, we could deduce the assertions corresponding to \cite[Lemmas 7 and 8]{MasYos05}:
\begin{itemize}
\item The flow $(X(t,\hat{v}),H(t,\hat{v}))_{t\in[0,t^{0}]}$ satisfies the condition $(\tilde{A}'-4)$ 
(as was seen in Section \ref{sec_mcm});

\item There exists a random variable $\Phi''_{t^{0}}\in L^{1}(\hat{P})$ such that
\begin{equation}
E\left[\sup_{|u|\ge B}\left|E[\hat{\psi}_{\ep}\exp(iuH_{t^{0}})|X_{0},X_{t^{0}}]\right|\right]
\le\frac{1}{B}\hat{E}[|\Phi''_{t^{0}}|]
\nonumber
\end{equation}
for every $B>0$.

\end{itemize}
The proof of Theorem \ref{th2} is thus complete.


\subsubsection*{Acknowledgement.} 
This work was partly supported by JSPS KAKENHI Grant Number 23740082 (H. Masuda).



\end{document}